\documentclass{article}

\usepackage[preprint]{neurips_2021}

\usepackage[utf8]{inputenc} 
\usepackage[T1]{fontenc}    
\usepackage{hyperref}       
\usepackage{url}            
\usepackage{booktabs}       
\usepackage{amsfonts}       
\usepackage{nicefrac}       
\usepackage{microtype}      
\usepackage{xcolor}         
\usepackage{amsmath}
\usepackage{dsfont}
\usepackage{graphicx}
\usepackage{array}
\usepackage{ulem}
\usepackage{wrapfig}
\usepackage{algorithm}
\usepackage{algorithmic}
\newcommand{\bX}{\boldsymbol{X}}
\newcommand{\bY}{\boldsymbol{Y}}
\newcommand{\bC}{\boldsymbol{C}}

\newcommand{\ba}{\boldsymbol{a}}
\newcommand{\bb}{\boldsymbol{b}}
\newcommand{\bc}{\boldsymbol{c}}
\newcommand{\bm}{\boldsymbol{m}}
\newcommand{\by}{\boldsymbol{y}}
\newcommand{\bx}{\boldsymbol{x}}
\newcommand{\bpi}{\boldsymbol{T}}

\newcommand{\bH}{\boldsymbol{H}}
\newcommand{\bI}{\boldsymbol{I}}

\newcommand{\bK}{\boldsymbol{K}}
\newcommand{\bp}{\boldsymbol{t}}
\newcommand{\bv}{\boldsymbol{v}}

\newcommand{\bu}{\boldsymbol{u}}
\newcommand{\setA}{\mathcal{A}}
\newcommand{\bB}{\boldsymbol{B}}

\newcommand{\bgamma}{\boldsymbol{\gamma}}

\def\bal#1{\begin{align}#1\end{align}}

\title{Unbalanced Optimal Transport through Non-negative Penalized Linear Regression}
\author{%
  Laetitia Chapel\thanks{First two authors have equal contribution} \\
  IRISA, Université Bretagne-Sud\\
  Vannes, France\\
  \texttt{laetitia.chapel@irisa.fr} \\
  \And
  Rémi Flamary\footnotemark[1] \\
  CMAP, Ecole Polytechnique\\
  Palaiseau, France\\
   \texttt{remi.flamary@polytechnique.edu}\\
   \And
   Haoran Wu \\
   LITIS \& IRISA\\
  Rouen \&  Vannes, France \\
  \texttt{haoran.wu@univ-ubs.fr} \\
   \And
   Cédric Févotte\\
   IRIT, Université de Toulouse, CNRS \\
   Toulouse, France \\
   \texttt{cedric.fevotte@irit.fr} \\
  \And
   Gilles Gasso \\
  LITIS, INSA Rouen Normandie \\
   Rouen, France\\
  \texttt{gilles.gasso@insa-rouen.fr} \\
}

\begin{document}

\maketitle

\begin{abstract}

This paper addresses the problem of Unbalanced Optimal Transport (UOT) in which
the marginal conditions are relaxed (using weighted penalties in lieu of
equality) and no additional regularization is enforced on the OT plan. In this
context, we show that the corresponding optimization problem can be reformulated
as a non-negative penalized linear regression problem. This reformulation allows
us to propose novel algorithms inspired from inverse problems and nonnegative
matrix factorization. In particular, we consider majorization-minimization which
leads in our setting to efficient multiplicative updates for a variety of
penalties. Furthermore, we derive for the first time an efficient algorithm to
compute the regularization path of UOT with quadratic penalties. The proposed
algorithm provides a continuity of piece-wise linear OT plans converging to the
solution of balanced OT (corresponding to infinite penalty weights). We perform
several numerical experiments on simulated and real data illustrating the new
algorithms, and provide a detailed discussion about more sophisticated
optimization tools that can further be used to solve OT problems thanks to our
reformulation. 

\end{abstract}


\section{Introduction}
\label{sec:itro}

Optimal Transport (OT) theory provides powerful tools for comparing probability
distributions and has been successfully employed in a wide range of machine
learning applications such as supervised learning \citep{frogner2015learning}, clustering \citep{ho2017}, generative modelling
\citep{arjovsky2017wasserstein}, domain adaptation \citep{courty2017optimal},
learning of structured data \citep{maretic2019got, titouan2019optimal} or natural language
processing \citep{kusner2015word}, among many others. {One 
reason
for those recent successes is the 
introduction of entropy-regularized OT that can be solved with the efficient
Sinkhorn-Knopp matrix scaling algorithm \citep{cuturi2013sinkhorn}.}
{
However, the classical OT problem seeks 
}the
optimal cost to transport \emph{all} the mass from a source distribution to a
target one \citep{villani03topics}, greatly limiting its use in scenarii where the
measures have different masses or when they contain noisy observations or
outliers. 

Unbalanced Optimal Transport (UOT) \citep{benamou2003numerical} has been introduced to tackle this
shortcoming, allowing some mass variation in the transportation problem. It is expressed as a relaxation of the Kantorovich formulation \citep{kantorovich1942}  by penalizing the
divergence between the marginals of the transportation plan and the given 
distributions.  
Several divergences can be considered, such as the Kullback-Leiber (KL)
divergence \citep{frogner2015learning,liero2018optimal}, the $\ell_1$ norm corresponding to the
partial optimal transport problem \citep{caffarelli2010, figalli2010}, or the
squared $\ell_2$ norm \citep{benamou2003numerical}. Regarding numerical solutions, \citet{chizat2018scaling} considered an entropic-regularized version of UOT leading to
a class of
scaling algorithms in the vein of the Sinkhorn-Knopp approach
\citep{sinkhorn1967concerning}. {The introduction of this entropic
regularization improves the scalability of OT, but involves a spreading of the
mass and a loss of sparsity in the OT plan. When a sparse transport plan is {sought}, the convergence is slowed down,
necessitating the use of acceleration strategies
\citep{thibault2021overrelaxed}}. 
Regarding UOT with the {(squared)} $\ell_2$ norm, \citet{blondel2018smooth} showed that the
resulting OT plan is sparse and proposed to use an
efficient L-BFGS-B algorithm \citep{Lbfgs_B_95} to address this case. Note that the L-BFGS-B method can be used to solve UOT with differentiable divergences even without the entropic-regularization on the OT
plan that induces the  Sinkhorn-like iterations.
Finally, also note that, as for 
balanced OT, UOT can 
be solved more efficiently when the data has a specific structure, {such as
{unidimensional} distributions \citep{Bonneel2019spot} or distributions supported on trees \citep{sato2020fast}.

\paragraph{Contributions.} In this paper, we show after some preliminaries that UOT can be
recast as a {convex} penalized linear regression problem with {non-negativity} constraints (Section \ref{sec:uot_linreg}).
The main interest of this reformulation resides in the fact that 
 non-negative linear regression 
has been
extensively studied in inverse problems and machine learning, offering a large panel of tools for devising new numerical algorithms. {Our reformulation involves a design/dictionary matrix that is structured and sparse. {Leveraging} this structure, we propose two new families of algorithms for solving the exact (i.e., without regularization of the plan) UOT problem in Section  \ref{sec:path}.}

{We first derive in Section \ref{sec:mmforuot} a new Majorization-Minimization (MM) algorithm for solving UOT with Bregman divergences, and more specifically KL and $\ell_2$-penalized UOT. The MM approach results in multiplicative updates that have appealing features: i) they are easy to
implement, ii) have low complexity per iteration and can be instantiated on GPU, iii) ensure monotonicity of the objective function and inherit existing convergence results. Our methodology is inspired by well-known algorithms in image restoration \citep{richardson1972bayesian,depi93} and non-negative matrix factorization (NMF) \citep{lee2001algorithm, dhillon2005generalized, fevotte2011algorithms}. Interestingly, the resulting multiplicative updates bear a similarity with the celebrated Sinkhorn scaling algorithm, with some key differences that are discussed.}


{Next, we derive in Section \ref{sec:regpath} an efficient algorithm to
compute the regularization path in $\ell_{2}$-penalized UOT. To do so, we build
on our proposed reformulation and more precisely on the fact that
$\ell_2$-penalized UOT can be reformulated as a weighted Lasso problem. We
propose a new methodology inspired by LARS
\citep{efron2004least,hastie2004entire}, which, to the best of our knowledge, is
the first regularization path algorithm for OT problems. It brings a novel
understanding of the properties of the evolution of the support of OT plans,
besides the practical interest of computing the complete regularization path when
hyperparameter validation is necessary.}


{Our new families of algorithms (MM for general UOT, LARS for
$\ell_{2}$-penalized UOT) are showcased in the numerical experiments of  Section
\ref{sec:expe}. Python implementation of the algorithms, {provided in supplementary}, will be released
with MIT license
on GitHub. The connection between UOT and linear regression that we
reveal in the paper opens the door to further fruitful developments and in
particular to more efficient algorithms, thanks to the large literature
dealing with non-negative penalized linear regression. We discuss those possible
research directions in Section \ref{sec:discuss}, before concluding the paper.}

\paragraph{Notations.} {Vectors such as  $\bm$  are written with lower case and bold
font, with coefficients $m_i$ or $[\bm]_i$, according to context. The $|\setA|$-dimensional sub-vector with indexes in set $\setA$ is written $\bm_{\setA}$. Matrices such as  $\boldsymbol{M}$ are written with upper case and bold
font, with coefficients $M_{i,j}$. We introduce a vectorization
operator defined by $\bm=\text{vec}(\boldsymbol{M})=[M_{1,1},M_{1,2},\dots, M_{n,m-1},
M_{n,m}]^\top$, i.e., the concatenation of the \emph{rows} of the
matrix, following the Numpy/C memory convention. 
$\mathds{1}_n$ is a
vector of $n$ ones and $\boldsymbol{M} \ge 0$ denotes entry-wise non-negativity. Finally, $D_\varphi$ is
the Bregman divergence generated by the strictly convex and differentiable
function $\varphi$, i.e., $D_\varphi(\bu,\bv)= \sum_{i} d_{\varphi}(u_{i},v_{i}) = \sum_i [
\varphi(u_i)-\varphi(v_i)-\varphi'(v_{i})(u_{i}-v_{i})]$.

}



\section{Reformulation of UOT as non-negative penalized linear regression}
\label{sec:uot_as_linreg}

\subsection{Background on Optimal Transport}
\label{sec:ot}

Let us consider two clouds of points $\bX=\{\bx_i\}^n_{i=1}$ and
$\bY=\{\by_j\}^m_{j=1}$. Let $\ba \in \mathbb{R}^{+}_{n}$ and $\bb
\in \mathbb{R}^{+}_{m}$ be two discrete distributions of mass on $\bX$ and $\bY$, such that $a_i$ (resp. $b_j$) is the mass at $\bx_i$ (resp. $\by_j$).
%
%
The \emph{balanced} OT problem, as defined by \citet{kantorovich1942}, is a linear
problem that computes the minimum cost of moving $\ba$ to $\bb$:
\begin{equation}
\text{OT}(\ba, \bb) = \min_{\bpi \geq 0} \langle \bC, \bpi\rangle\quad 
\text{such that (s.t.)} \quad \bpi \mathds{1}_m = \ba, \bpi^\top \mathds{1}_n = \bb
\label{eq:exactot}
\end{equation}
where $ \langle \cdot, \cdot \rangle$ is the Frobenius inner product, $\bpi \in \mathbb{R}^{+}_{n\times m}$ is the \emph{transport plan} and $\bC \in \mathbb{R}^{+}_{n\times m}$ is the \emph{cost matrix}. The entry $C_{i,j}$ of $\bC$ represents the cost of moving point $\bx_i$ to $\by_j$. The Wasserstein 1-distance (also known as the earth
mover's distance) is obtained for $C_{i,j}=\|\bx_i-\by_j\|$.
The constraints on the
transport plan $\bpi$ require that $\|\ba\|_1 = \|\bb\|_1$ and that
\textit{all} the mass from $\ba$ is transported to $\bb$.
These constraints can be alleviated through relaxation,
leading to UOT  \citep{benamou2003numerical}:
\begin{eqnarray}
\text{UOT}^{\boldsymbol{\lambda}}(\ba, \bb) = \min_{\bpi \geq 0} \quad \langle \bC, \bpi\rangle
+ \lambda_1 D_\varphi(\bpi \mathds{1}_m, \ba) +  \lambda_2 D_\varphi(\bpi^\top
\mathds{1}_n, \bb).
\label{eq:UOTdiv}
\end{eqnarray}
The deviations from the true marginals are penalized by means of a given Bregman divergence $D_\varphi$, as introduced in \citet{chizat2018scaling},
where $\lambda_1$ and $\lambda_2$ are hyperparameters that represent the strengths
of penalization. Note that balanced OT \eqref{eq:exactot} is recovered when $\lambda_1 = \lambda_2 \to \infty$. Furthermore, when $\lambda_1$ \emph{or} $\lambda_2 \to \infty$, we recover semi-relaxed OT \citep{rabin2014adaptive}. In practice, authors often set $\lambda_1 = \lambda_2= \lambda$ for UOT in order to reduce the necessity of hyperparameter tuning. Various divergences have been considered in the literature. The
$\ell_1$ norm gives rise to so-called \textit{partial} optimal transport
\citep{caffarelli2010}. The squared $\ell_2$ norm provides a sparse and smooth
transport plan \citep{blondel2018smooth} when introducing a strongly convex
term in Eq.~\eqref{eq:UOTdiv}. \citet{chizat2018scaling} derive efficient
algorithms to solve Eq.~\eqref{eq:UOTdiv} for several divergences by adding an additional
regularization term $\lambda_{\text{reg}} D_\varphi(\bpi,\ba \bb^\top)$. In particular, entropic regularization is obtained when the KL divergence is used, promoting a dense transport plan unlike exact UOT.


\subsection{Reformulation of UOT}
\label{sec:uot_linreg}

\paragraph{UOT cast as regression.} Let $\bp=\text{vec}(\bpi)$, $\bc=\text{vec}(\bC)$ and  $\by^\top=[\ba^\top,\bb^\top]$. Problem~\eqref{eq:UOTdiv} can be re-written as
        \begin{equation}
            \min_{\bp\geq 0}\quad  F_\lambda(\bp) \stackrel{\text{def}}{=} \frac{1}{\lambda}\bc^\top\bp + D_\varphi(\bH\bp,\by) \label{eq:linreg}
            \end{equation}
and as such be expressed as a non-negative penalized linear regression
problem, where the {\em design matrix} $\bH=[\bH_r^\top,\bH_c^\top]^\top$ is the
concatenation of the matrices $\bH_r$ and $\bH_c$ that compute sums of the rows
and columns of $\bpi$, respectively (see expressions in Section~\ref{sec:design}
of the supplementary material). Note that, for the sake of simplicity, we
consider here $\lambda_1 = \lambda_2 = \lambda$ but this hypothesis could be
easily alleviated for a given family of divergences (see Sec.~\ref{sec:discuss}
for a discussion).
%
%
%
Important features of  Eq.~\eqref{eq:linreg} should be discussed. {First, $F_\lambda(\bp)$ is convex thanks to the convexity of Bregman divergences w.r.t. their first argument.} Second, $\bH$ is very structured and sparse (with a ratio of only
$\frac{1}{m+n}$ non-zero coefficients) which will allow {for more efficient
computations and updates than with a dense $\bH$.}
 Finally, since $\bp\geq 0$ and $\bc\geq 0$, the linear term can be expressed as $\frac{1}{\lambda}\bc^\top\bp=\frac{1}{\lambda}\sum_i
c_i t_i=\frac{1}{\lambda}\sum_i
c_i |t_i|$. This corresponds to a weighted $\ell_1$ regularization, promoting sparsity in $\bp$ and hence in the transport plans. 
Note that the ``sparse'' regularization is here controlled by $\frac{1}{\lambda}$ (instead of $\lambda$ in
classical penalized linear regression), meaning that the sparsity promoting term will be more aggressive for small $\lambda$. 




\paragraph{Solving problem~\eqref{eq:linreg}.} 
{Problems of the form of Eq.~\eqref{eq:linreg} are well-known in inverse
problems and NMF. In inverse problems, $\bp$ typically acts as a clean image
degraded by operator $\bH$ (e.g., a convolution) and noise. The data fitting
term $D_\varphi(\bH\bp,\by)$ captures assumptions about the noise corrupting the
observed image $\by$. Sparsity is a common regularizer of $\bp$. In NMF, given a
set of nonnegative samples $\{ \by_{l} \}$ one wants to learn a non-negative
dictionary $\bH$ and non-negative lower-dimensional embeddings $\{ \bp_{l} \}$
such that $\by_{l} \approx \bH \bp_{l}$ \citep{lee1999learning}. Updating the
latter involves optimization problems of form~\eqref{eq:linreg} (with or without
sparse regularization). In contrast to problem \eqref{eq:linreg}, the data
fitting term is more commonly $D_\varphi(\by, \bH\bp)$ instead of
$D_\varphi(\bH\bp,\by)$ in inverse problems and NMF. This is because the former
is a log-likelihood in disguise for the mean-parametrized exponential family,
and takes important noise models as special cases, such as Poisson, additive
Gaussian or multiplicative Gamma noise \citep{fevotte2011algorithms}. Using such
penalizations with reversed arguments would be possible in our case as well but
we stick to the now standard formulation of \citep{liero2018optimal, chizat2018scaling}
for simplicity.

In the next section, we will first leverage a classical family of algorithms in inverse problems and NMF, namely MM, to obtain new algorithms for KL and $\ell_{2}$-penalized UOT (possibly with entropic regularization in the first case). Second, we will leverage results about non-negative Lasso  to design an efficient algorithm to compute the regularization path of $\ell_{2}$-penalized UOT.}

\section{{Novel} numerical solvers for UOT}
\label{sec:path}

\subsection{Majorization-Minimization (MM) for UOT}
\label{sec:mmforuot}

{
\paragraph{General MM framework.}

MM algorithms have been around a long time in inverse problems and NMF to solve problems of form~\eqref{eq:linreg}. Classical algorithms for NMF such as \citep{lee2001algorithm} have built on seminal MM algorithms for inverse problems such as \citep{richardson1972bayesian,depi93}. Subsequent works in NMF such as \citep{dhillon2005generalized, fevotte2011algorithms, yang11} have further contributed novel MM algorithms for larger classes of problems, including larger families of divergences. In a nutshell, MM consists in iteratively building and minimizing an upper bound of the objective function which is tight at the current parameter estimate (and referred to as \emph{auxiliary function}), see \citet{hunter2004tutorial,Sun2017} for tutorials. In NMF, a common approach consists of alternating the updates of the dictionary $\bH$ and of the embeddings. In our case, $\bH$ is fixed and we may use the results of \citep{dhillon2005generalized} to build an auxiliary function for term $D_\varphi(\bH\bp,\by)$, to which we may simply add the linear term $\bc^{\top} \bp/\lambda$ to obtain a valid auxiliary function for $F_\lambda (\bp)$. Let $\tilde{\bp}$ denote the current estimate of $\bp$, $\tilde{Z}_{i,j}=\frac{H_{i,j}\tilde{t}_j}{\sum_l H_{i,l}\tilde{t}_l}$ and
\begin{equation}
	G_\lambda (\bp, \tilde{\bp}) = \sum_{i,j} \tilde{Z}_{i,j} \varphi\left(\frac{H_{i,j}t_j}{\tilde{Z}_{i,j}}\right) + \sum_j \left[ \frac{c_{j}}{\lambda}  - \sum_{i} H_{i,j} \varphi'(y_{i}) \right] t_{j}	+ cst,
	\label{eq:auxiliary_func}
\end{equation}
where $cst = \sum_{i}[\varphi'(y_{i})y_{i} - \varphi(y_{i}] $. Then, $G_\lambda (\bp,\tilde \bp)$ is an auxiliary function for $F_\lambda (\bp)$, i.e., $\forall \bp$, $G_\lambda (\bp,\tilde{\bp}) \ge F_\lambda (\bp)$ and $G_\lambda (\tilde{\bp},\tilde{\bp}) = F_\lambda (\tilde{\bp})$. Let $\bp^{(k+1)} = \text{argmin}_{\bp \ge 0} G_\lambda ({\bp}, \bp^{(k)} )$, then $ F_\lambda (\bp^{(k)}) =  G_\lambda (\bp^{(k)}, \bp^{(k)} ) \ge G_\lambda (\bp^{(k+1)}, \bp^{(k)} ) \ge F_\lambda (\bp^{(k+1)})$, producing a descent algorithm over $F$. The trick to obtain $G$ is to apply Jensen inequality to $\varphi( \sum_{j} H_{i,j} t_{j}) =\varphi( \sum_{j} \tilde{Z}_{i,j} \frac{H_{i,j}}{\tilde{Z}_{i,j}} t_{j}) \le \sum_{j} \tilde{Z}_{i,j} \varphi(\frac{H_{i,j}}{\tilde{Z}_{i,j}} t_{j}) $, thanks to the convexity of $\varphi$, see details in \citep{dhillon2005generalized}. We provide below the resulting algorithms for the KL and $\ell_{2}$ penalizations, with
detailed computations available in Section \ref{sec:app_mult} of the supplementary.

\paragraph{MM for KL-penalized UOT.} The KL divergence is obtained with $\varphi (y) = y \log y - y$. Minimizing $G_\lambda (\bp,\bp^{(k)})$ in that case leads to following multiplicative update:
 \begin{equation}
 	t_j^{(k+1)} = t^{(k)}_j \, \exp\left(\frac{[\bH^\top \log(\by) - \bH^\top \log\left(  \bH \bp^{(k)} \right)]_j - \frac{1}{\lambda} c_{j}}{[\bH^\top \mathds{1} ]_j}\right). 
\end{equation}
Owing to the structure of $\bp$ and $\bH$, the update can be re-written in the following matrix form:
%
\begin{equation}
	\bpi^{(k+1)}= \text{diag}\left(\frac{\ba}{\bpi^{(k)} \mathds{1}_m}\right)^{\frac{1}{2}}\left(\bpi^{(k)} \odot \text{exp}\left(-\frac{\bC}{2\lambda}\right)\right) \text{diag}\left(\frac{\bb}{{\bpi^{(k)\top}}\mathds{1}_n}\right)^{\frac{1}{2}},
	\label{eq:mult_update_kl}
\end{equation}
where $\odot$ is entrywise multiplication and divisions are taken entrywise as well. The multiplicative update \eqref{eq:mult_update_kl} is {remarkably} 
similar
to the well-known Sinkhorn-Knopp algorithm 
that has
been used in numerous OT problems involving KL regularization. 
But {instead of two separate steps for the left and right scaling, Eq.~\eqref{eq:mult_update_kl} applies these scalings simultaneously in a unique update using 
the diagonal matrices (and a form of geometrical average).} Also note how the scaling
factor $\text{exp}\left(-\frac{\bC}{2\lambda}\right)$ penalizes along 
iterations the coefficients of the transport plan with large costs. 

\paragraph{MM for $\ell_2$-penalized UOT.} The quadratic loss is obtained with $\varphi(y) = \frac{y^2}{2}$. In that case, minimizing $G_\lambda (\bp,\bp^{(k)})$ s.t. non-negativity leads to following multiplicative update:
\begin{equation}
	\bpi^{(k+1)}= \bpi^{(k)} \odot \frac{\max\left(0,  \ba \mathds{1}_m^\top + \mathds{1}_n \bb^\top  -\frac{1}{\lambda} \bC \right)}{\bpi^{(k)} \mathbf{O}_m + \mathbf{O}_n \bpi^{(k)}} \quad \text{with} \quad \mathbf{O}_\ell = \mathds{1}_\ell \mathds{1}_\ell^\top.
	\label{eq:mult_update_quadratic_divergence_l2}
\end{equation}
{Interestingly enough, update \eqref{eq:mult_update_quadratic_divergence_l2} prunes any coefficient $T_{i,j}$ in $\bpi$ such that $a_i+b_j-\frac{1}{\lambda}C_{i,j}<0$ from the very first iteration, providing a useful certificate on the support of the solution.


\subsection{Regularization path for $\ell_2$-penalized UOT}
\label{sec:regpath}


{Let us focus on the case where  $D_{\varphi}$ is a quadratic divergence. As
mentioned  in Section \ref{sec:uot_linreg}, Eq.~\eqref{eq:linreg} is then a
positive weighted Lasso problem, allowing us to derive
the first regularization path algorithm for computing the whole set of solutions
for a varying $\lambda$ from 0
to $+\infty$. Note that the path's extreme point recovers the balanced OT solution.
We show
that the path is piecewise linear in $1/\lambda$ between changes in the
active set {$\mathcal{A}=\text{supp}(\bp^{\lambda})$, where $\bp^{\lambda} = \text{vec}(\bpi^{\lambda})$ and $\bpi^{\lambda}$ is the OT plan for given hyperparameter $\lambda$.}
The main steps of the algorithm are roughly as follows: given a current solution
$(\lambda_k, \bpi^{\lambda_k})$ and a current active set $\mathcal{A}_k$, we
look for the next value $\lambda_{k+1}> \lambda_k$ such that the active set
changes (i.e., $\mathcal{A}_{k+1} \neq \mathcal{A}_k$), either because one component
enters or leaves the active set. 
We describe our algorithm below.}

\paragraph{KKT conditions of the $\ell_2$-penalized UOT problem.}
The Lagrangian for problem \eqref{eq:linreg} writes:
\begin{equation}
	{L}_\lambda(\bp, \bgamma) = \frac{1}{\lambda}\bc^\top\bp + {\frac{1}{2}}(\bH\bp-\by)^\top(\bH\bp-\by)-\bgamma^\top\bp
\end{equation}
where $\bgamma$ represents the Lagrange parameters. We denote 
$\bm =  \bH^\top\by= \text{vec}(\ba {\mathds{1}_m^\top}+ {\mathds{1}_n} \bb^\top)$. KKT optimality conditions state that i) $\nabla_{\bp} {L}_\lambda(\bp, \lambda) =\frac{1}{\lambda}\bc + \bH^\top\bH\bp-\bm - \bgamma = 0$ (stationarity condition), ii) $\bgamma \odot \bp = 0$ {(complementary condition)} and iii) $\bgamma \geq 0$ (feasibility condition). 

\paragraph{Piecewise linearity of the path.} Assume that, at iteration $k$, we
know the current active set $\setA = \mathcal{A}_k$ and we look for
$\bp_{\setA}^{\lambda} $ (the other values of $\bp_{\setA}$ being $0$). Let 
$\bH_{\setA}$, $\bm_A$ and $\bc_{\setA}$ denote the corresponding sub-matrix and vectors (see
Appendix \ref{sec:matrices_reg_path} for rigorous definitions). Because of the
complementary condition, we have $\bgamma_{\setA} = \boldsymbol{0}$. {Using $\lambda =
\lambda_k+\epsilon$, with $\epsilon > 0$ small enough to ensure that the
active set remains the same, the stationarity condition writes:}
\begin{eqnarray}
 \bH_{\setA}^\top\bH_{\setA}\bp^{\lambda}_{\setA} = \bm_{\setA} - \frac{1}{\lambda}\bc_{\setA}  & \Rightarrow & \bp_{\mathcal{A}}^{\lambda} = \tilde{\bm}_{\setA} - \frac{1}{\lambda}\tilde{\bc}_{\setA}
	\label{eq:opti}
\end{eqnarray}
with $\tilde{\bm}_{\setA} =  (\bH_{\setA}^\top\bH_{\setA})^{-1}{\bm}_{\setA}$ and $\tilde{\bc}_{\setA} =  (\bH_{\setA}^\top\bH_{\setA})^{-1}{\bc}_{\setA}$ . {Eq.~\eqref{eq:opti} shows that the optimal $\bp_{\mathcal{A}}^{\lambda}$ (and hence $\bp^{\lambda}$) can be solved for any $\lambda \in [\lambda_{k}, \lambda_{k+1}]$, i.e., when the active set $\mathcal{A}$ remains the same, by solving a linear problem. It also reveals the piecewise linearity in $\lambda^{-1}$ of the path when $\mathcal{A}$ is fixed. As expected, balanced OT is recovered when $\lambda \to \infty$.}

\paragraph{Finding $(\lambda_{k+1},\mathcal{A}_{k+1}) $ given $(\lambda_k,\mathcal{A}_k)$.}
Given a current solution $(\lambda_k, \bp^{\lambda_k})$ and $\lambda = \lambda_k+\epsilon$, we
increase the $\epsilon$ until we reach a change in the set of active components. This happens whenever the first of the following two situations occurs.

$\bullet$ \textbf{One component in $\mathcal{A}$ becomes inactive.} 
In that case, we remove the index $i \in \mathcal{A}$ with the smallest $\lambda_r>\lambda_k$ that violates the constraint. In such case, {$[\tilde{\bm}_{\setA}]_{i} =  [\tilde{\bc}_{\setA}]_{i}/\lambda$} and we may write
\begin{eqnarray}
  \lambda_r = \min_{>\lambda_k} \left( \frac{\tilde{\bc}_{\setA}}{\tilde{\bm}_{\setA}  } \right)
  \label{eq:lambdar}
\end{eqnarray}
where  $\min_{>\lambda_k}$ indicates the minimum value in the vector greater than $\lambda_k$ and the division is entrywise.

$\bullet$ \textbf{One component in $\bar{\mathcal{A}}$ becomes active.} This occurs when the KKT positivity constraint $\bgamma_{\bar{\setA}} \geq \mathbf{0}$ becomes violated. Assume this happens at index $i \in \bar{\mathcal{A}}$ for the smallest value $\lambda_a>\lambda_k$ of $\lambda$. {In such case, the stationarity condition outside the active set can be rewritten}:
\begin{equation}
\left[\frac{1}{\lambda}\bc_{\bar{\setA}}+  \Big[\bH^\top \bH\big(\tilde{\bm} + \frac{1}{\lambda}\tilde{\bc}\big)\Big]_{\bar{\setA}}-    \bm_{\bar{\setA}}\right]_{i} = [\bgamma_{\bar{\setA}}]_{i} 
\ \Rightarrow \ 
\lambda_a =  \min_{>\lambda_k}\left(\dfrac{\bc_{\bar{\setA}} - \big[\bH^\top\bH\tilde{\bc}\big]_{\bar{\setA}} }{\bm_{\bar{\setA}}-\big[\bH^\top\bH\tilde{\bm}\big]_{\bar{\setA}} }\right),
\label{eq:lambdaa}
\end{equation}
where $\tilde{\bm}$ (resp. $\tilde{\bc}$)  equals $\tilde{\bm}_{\setA}$ (resp. $\tilde{\bc}_{\setA}$) on $\setA$ and zero on $\bar{\setA}$.

In practice, at each { step of the path}, we compute both $\lambda_r$ and $\lambda_a$, set $\lambda_{k+1}= \min\{ \lambda_r,\lambda_a \}$ and update the active set  accordingly. 
\paragraph{Numerical computation of the entire path.}
Eq.~\eqref{eq:opti} involves the computation of the matrix $(\bH_{\setA}^\top\bH_{\setA})^{-1}$, which is of size $|\setA|\times|\setA|$. As only one index leaves or enters the active set at each iteration, we can use the Schur complement of the matrix to compute its value from $(\bH_{\setA_{k}}^\top\bH_{\setA_k})^{-1}$, alleviating the computational burden of the algorithm as it only involves matrix-vector computations (see Section \ref{sec:matrices_reg_path} of supplementary).
Algorithm \ref{alg:regpath} sums up the different steps of the full path computation.
At each iteration, we compute $\lambda_a$,  $\lambda_r$, update the inverse matrix $(\bH_{\setA_{k}}^\top\bH_{\setA_k})^{-1}$ and estimate the solution
$\bp^{\lambda_{k+1}}$ with a complexity of $O(n m)$. 

\begin{algorithm}[t]
	\caption{Regularization path of $\ell_2$-penalized UOT}
	\begin{algorithmic} 
	\STATE Require: $\ba$, $\bb$, $\bC$, $\lambda_0=0$, $\bp_0 = \boldsymbol{0}$, $\setA = \setA_0=\emptyset$, $k=1$
	\STATE $\lambda_1=\min \dfrac{\bc_{\bar{\setA}}}{\bm_{\bar{\setA}}}$, $\setA = \setA_1 = \arg\min \dfrac{\bc_{\bar{\setA}}}{\bm_{\bar{\setA}}}$, $\bH_{\setA}^\top\bH_{\setA}=2$
	\STATE $\bp^{\lambda_1}_{\setA_1 }=\frac{\bm_{{\setA}}}{2}-\frac{1}{\lambda_1}\frac{\bc_{\setA}}{2}$
	\WHILE{$(\bH\bp^{\lambda_k} - \by)^\top(\bH\bp^{\lambda_k} - \by)\neq 0$}
	\STATE $\lambda_r$, $\lambda_a$ $\leftarrow$ Compute  as in Eq. (\ref{eq:lambdar}) and  Eq. (\ref{eq:lambdaa})
	\STATE $\lambda_{k+1} \leftarrow \min(\lambda_r, \lambda_a)$
	\STATE $\bp^{\lambda_{k+1}}_{\setA} \leftarrow (\bH_{\setA}^\top\bH_{\setA})^{-1}\bm_{\setA}-\frac{1}{\lambda_{k+1}}(\bH_{\setA}^\top\bH_{\setA})^{-1}\bc_{\setA}$
	\STATE $\setA = \setA_{k+1} \leftarrow$ Update active set for next iteration.
	\STATE  $(\bH_{\setA}^\top\bH_{\setA})^{-1}\leftarrow$ Update from
	$(\bH_{\setA_k}^\top\bH_{\setA_k})^{-1}$ with Schur complement (see supplementary \ref{sec:matrices_reg_path})
	\STATE $k \leftarrow k + 1$
	\ENDWHILE
	\RETURN $(\lambda_k, \bp^{\lambda_k})_k$
	\end{algorithmic}
	\label{alg:regpath}
	\end{algorithm}
	
\paragraph{Regularization path of the semi-relaxed $\ell_2$-penalized UOT.}
As a side result, let us consider the semi-relaxed OT problem $
		\text{SROT}^{{\lambda}}(\ba, \bb) = \min_{\bpi \geq 0,  \bpi^\top\mathds{1}_n\ = \bb} \,\, \langle \bC, \bpi\rangle
		+ \lambda \|\bpi \mathds{1}_m -  \ba\|^2$. 
		The main difference with UOT is that the equality constraint $\bpi^\top\mathds{1}_n\ =\bb$ (equivalent to $\bH_c\bp = \bb$) must always be met. This leads to the following Lagrangian:
\begin{equation}
	L_\lambda(\bp,\bgamma,\bu) = \frac{1}{\lambda}\bc^\top\bp + {\frac{1}{2}}(\bH_r\bp-\ba)^\top(\bH_r\bp-\ba)+ (\bH_c\bp-\bb)^\top\bu-\bgamma^\top\bp,
	\label{eq:lagrang_SROT}
\end{equation}
where $\bu \in \mathbb{R}^m$ contains {the Lagrange parameters associated to the $m$
equality constraints}. The KKT optimality conditions now dictate that i)
$\nabla_{\bp} {L}_\lambda(\bp,\bgamma,\bu) =\frac{1}{\lambda}\bc +
\bH_r^\top\bH_r\bp-\bH_r^\top\ba +\bH_c^\top\bu- \bgamma = 0$, ii) $\bgamma \odot \bp = 0$, iii) $\bgamma  \geq 0$ and $\bH_c\bp-\bb=\boldsymbol{0}$. 
We can use the same reasoning than previously
to compute the entire path. Details are provided in Section
\ref{sec:reg_path_semi} of the supplementary. {The main difference lies in solving, at each
iteration, a linear system of size $(m + |\setA|)$ to comply with  the marginal
equality constraint. The path is initialized as follows: the $j^{th}$ column of $\bpi^0$ for $\lambda_0=0$ is set to the weighted canonical vector $b_{i^{\star}}\mathbf{e}_{i^{\star}}$, where $i^{\star} = \text{argmin} \{C_{i,j}\}_{i} $.

\section{Numerical experiments}
\label{sec:expe}

In this section, we first 
 show the solutions obtained with our solvers on simple and
interpretable examples. We then evaluate the computational complexity of the
different algorithms and finally we show an application where the regularization path can be
used on a domain adaptation problem.

\paragraph{Illustration of the algorithms.}
\label{sec:}
We first illustrate the regularization path for $\ell_2$-penalized UOT on a simple example between two distributions containing 3
points each, with different masses and a cost matrix $\bC$ given in Fig.
\ref{fig:path_simple} (left). We can see on Fig.
\ref{fig:path_simple} (right) that, starting from $\lambda_{0}=0$ and $\bpi = 0$,  we  successively add or remove components in the active
set $\setA$ when increasing the $\lambda$ values. When $\lambda=\infty$, we
recover the balanced OT solution. Recall that the path is linear in $1/\lambda$ (and not $\lambda$).
We then illustrate the  path for both  $\ell_2$-penalized UOT and semi-relaxed
UOT on two 2D distributions with $n=m=100$ samples. We can see in Fig.~\ref{fig:path_2d} 
the difference between the two regularization paths for specific values of $\lambda$. UOT starts with an empty plan for $\lambda=0$ and then activates
samples from both source and target from the closest to the farthest ones until
convergence to the balanced OT plan. Semi-relaxed UOT starts with all
target samples active due to marginal constraints and progressively activates the source samples.

\begin{figure}[t]
    \centering
\includegraphics[width=1\textwidth]{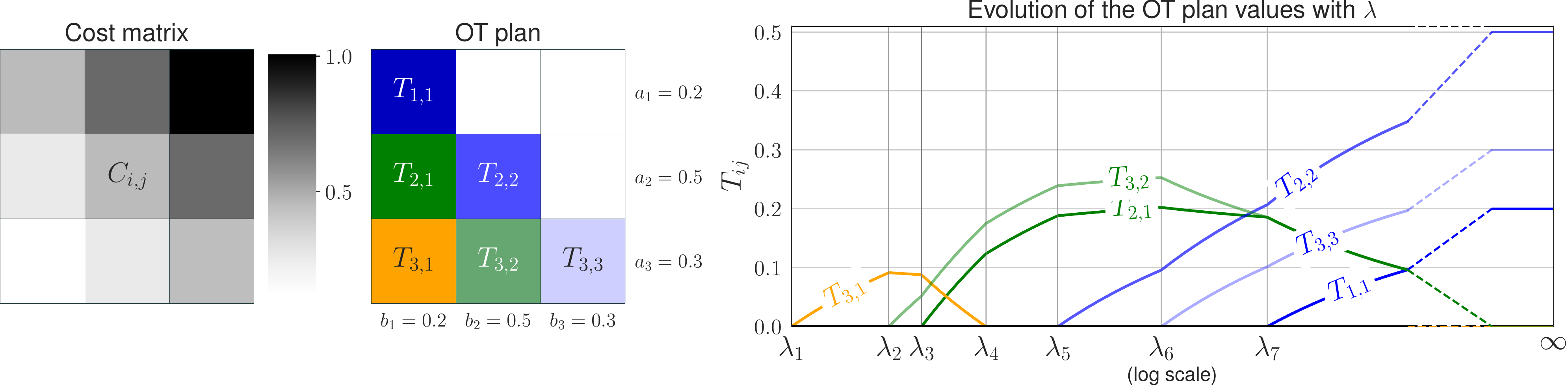}
\vspace*{-4mm}
    \caption{(Left) cost matrix $\bC$ (the higher the cost, the darker the color); (middle) OT plan whose cells are color-coded with respect to the $\lambda$ values at which they are activated. The blank cells never enter the active set as the corresponding cost it too high; (right) evolution of $T_{i,j}$ when $\lambda$ increases. Note that the $x$-axis is in log scale and is discontinued between $\lambda_7$ and $\infty$.}
    \label{fig:path_simple}
\end{figure}

\begin{figure}[t]
    \centering
    \includegraphics[width=0.99\textwidth]{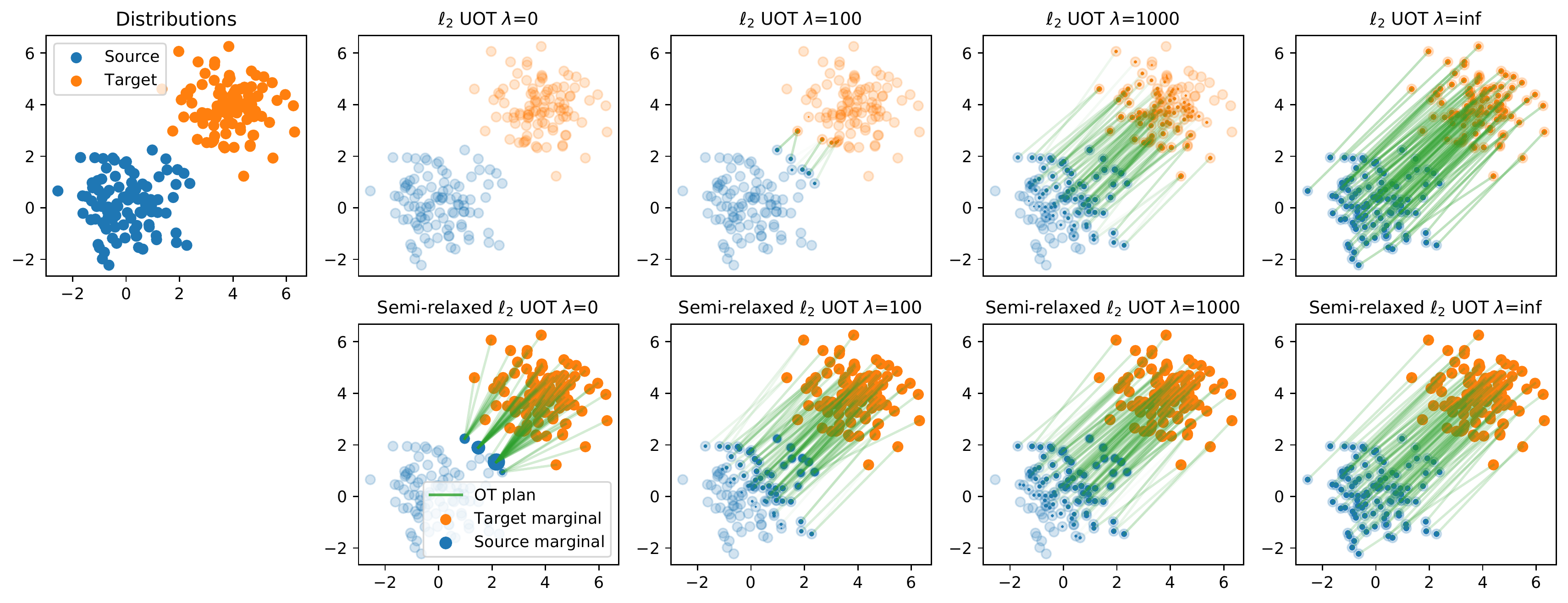}
    \caption{Regularization paths for 2D empirical distributions for $\ell_2$-penalized UOT
    (top) and semi-relaxed UOT (bottom). The OT plan is shown as green lines
    between the source and target samples when $T_{i,j}>0$ and the
    resulting marginals are shown as filled circles.}
    \label{fig:path_2d}
\end{figure}
}

\paragraph{Comparison of the performances of the algorithms.}
We now provide an empirical evaluation of the running times of the proposed
algorithms, using 2 sets of 10-dimensional points with $n=m$ and drawn according
to IID Gaussian distributions.  The cost matrix $\bC$ is computed using a squared $\ell_2$
norm. We first study the running times of the regularization path algorithm, for $n=m$
ranging from 100 to 1000, averaging the results over 5 runs, see Fig.
\ref{fig:timings} (left). We empirically observe that log-log plot is near-linear, with an empirical complexity $O(n^{3.27})$ in this example.

Using  $n=m=500$, we compare the running times of the current state-of-the-art BFGS algorithm
\citep{blondel2018smooth} using SciPy \citep{2020SciPy-NMeth} and those of our
algorithms: the $\ell_2$-penalized UOT formulated as a Lasso problem (with both
the Celer algorithm \citep{pmlr-v80-massias18a} and the coordinate descent solvers from
Scikit-learn \citep{pedregosa2011scikit}), the
multiplicative algorithm for both the $\ell_2$ and the KL penalties and the
regularization path algorithm (see Section \ref{sec:details_expes} of the
supplemental material
for more details about the solvers and their parameters). Figure
\ref{fig:timings} (middle and right) shows the average running time for all
algorithms. For $\ell_2$-penalized UOT, we observe that, for large $\lambda$ values,
the Lasso solvers are the fastest and that, whatever the value $\lambda$, BFGS is the slowest. We also notice that, for large $\lambda$, the running times for computing the path remain constant: when the last active set is found,
computing the OT plan only involves a weighted sum. As for KL-penalized UOT, the BFGS algorithm is more efficient when large values of $\lambda$ are considered. One can also notice that, 
similarly to Sinkhorn which is fast for large regularization values, the multiplicative algorithms for both penalties are also fast for high $1/\lambda$ values. 
\begin{figure}[!t]
    \includegraphics[width=1\textwidth]{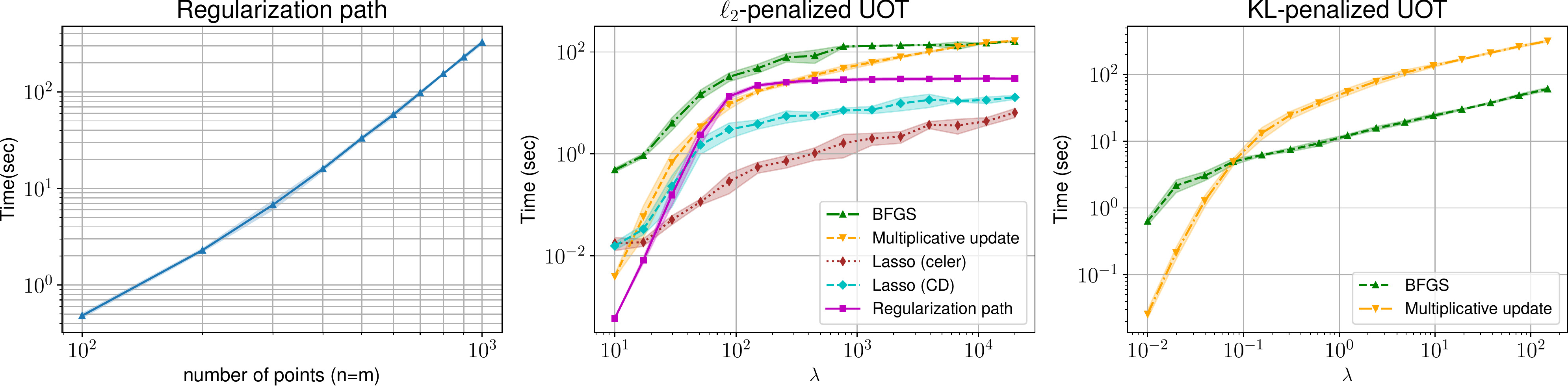} 
    \vspace*{-4mm}
    \caption{(Left) Running times of Alg.~1 w.r.t. the number of points; (middle) comparison of $\ell_2$-penalized UOT with $m=n=500$ (right) likewise for KL-penalized UOT. Dark curves (resp. shaded regions) represent average (resp. variance) values over 5 runs. }
    \label{fig:timings}
\end{figure}

\paragraph{Regularization path for unbalanced domain adaptation.}

We demonstrate the interest of having the entire regularization path in a
classification context where some of the data collection may be polluted by
outliers. We consider a setup similar to \citet{mukherjee2020outlier}. Let the
source $\bX$ be a set of 400 MNIST digits sampled from the digits $0,1,2,3$
(100 points per class) and let the target $\bY$ be a set of digits $0,1$ of
MNIST~\citep{lecun2010mnist} and of digits $8,9$ from Fashion MNIST~\citep{xiao2017fashion}. Our setting is  
simple classification: we classify a sample of the target dataset by
propagating \begin{wrapfigure}{r}{0.35\textwidth}\vspace*{-2mm}
  \begin{center}
    \includegraphics[width=0.3\textwidth]{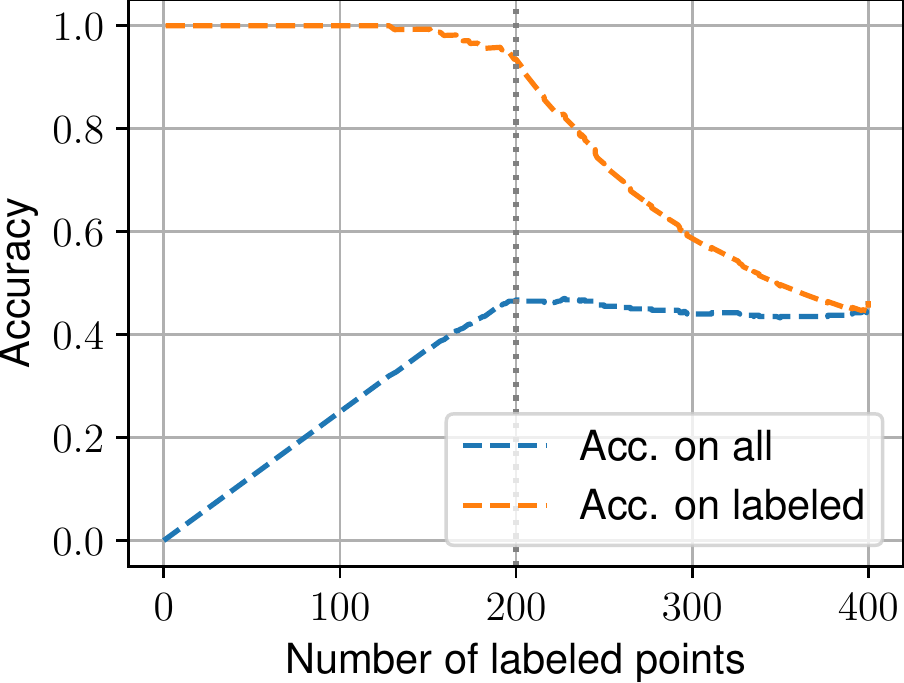}
  \end{center}
      \vspace*{-4mm}
  \caption{Evolution of the classification accuracy for the domain adaptation problem w.r.t. the number of classified points.}\vspace*{-4mm}
  \label{fig:classifMNIST}
\end{wrapfigure} the label of the source sample it is the most transported to, provided that
the transported mass of the target point is greater than $0.25b_j$. {Note that similarly to \citet{mukherjee2020outlier} a validation set can be used
here to select the best $\lambda$.}
Figure
\ref{fig:classifMNIST} shows the overall accuracy, defined as the number of
samples that are correctly classified divided by the total number of points, and
the current accuracy, which is the proportion of well-classified points among
the points that are classified, i.e., that are receiving mass. One can notice that, as the number of classified
points increases (with $\lambda$), the overall accuracy increases as more and
more points are well classified while the current accuracy remains stable until
outliers are included in the labeled set. {This suggests that UOT can be used
not only for classification but also as an automated outlier detection method.
}


%

\section{Discussion and perspectives}
\label{sec:discuss}

We showed that UOT can be recast as a non-negative penalized linear
regression problem, encouraging us to dig into this {well-established} field of research in order to {adapt existing algorithmic solutions to the structure of the UOT problem}. 
In this section, we discuss the relation between the proposed
algorithms and classical solvers used in OT, and also investigate 
{some} research
directions that can {widen the scope of proposed methods.} 

{
\paragraph{Multiplicative algorithms for UOT.} As discussed in Section
\ref{sec:mmforuot}, the multiplicative updates for the KL divergence obtained from MM resemble the
Sinkhorn algorithm from \citet{chizat2018scaling}, except for the joint scaling and the weighting matrix $\exp(-\bC/2)$. Interestingly, this
scaling matrix also appears in the Inexact Proximal Point OT (IPOT) algorithm of \citet{xie2020fast} to solve balanced OT. As a matter of fact, we show in Section~\ref{sec:prox_mm} of the supplementary that IPOT is a MM algorithm. The idea is to re-write the OT objective as $[\langle \bC, \bpi\rangle
+\lambda D_\varphi(\bpi,\ba\bb^\top)] -\lambda D_\varphi(\bpi,\ba\bb^\top)$ and upper bound the concave term by its tangent. This further supports the interest of MM for OT and UOT, and highlight an important feature of one of our contributions: designing the first Sinkhorn-like multiplicative algorithm for UOT that can be applied when the OT plan is not entropy-regularized.



\paragraph{More efficient solvers.} Despite the
positive experimental results of Section~\ref{sec:expe}, multiplicative and regularization path algorithms can be slow, especially for large
values of $\lambda$. Various accelerations can be envisaged. Regarding path algorithms, the approach of \citet{mairal2012complexity} can compute a regularization path with precision $\epsilon$ in $o(1/\epsilon)$ iterations. This would lead in our setting to a full complexity of $O({mn}/{\epsilon})$ that is even interesting to approximate balanced OT. Another way to speed up computations is to use \emph{screening}. In sparse regression, this consists of eliminating during optimization components that will not belong to the support of the solutions thanks to safe screening tests. Methods such as \citep{ghaoui2010safe, wang2015lasso, Dantas_2021_safe_screening_KLDiv} can readily be adapted to our $\ell_2$ or KL-penalized UOT algorithms. Finally, an other line of improvement is to consider stochastic optimization methods such as \citep{defazio2014saga}. Given the particular structure of $\bH$, the complexity of stochastic updates shall be small and can lead to very efficient implementations \citep{nesterov2014subgradient}.

\paragraph{General case and entropy-regularized UOT.}
{Following \citep{frogner2015learning,chizat2018scaling,sejourne2019sinkhorn}, general regularized UOT  can be expressed as:}
\begin{eqnarray}
    \text{RUOT}^{\boldsymbol{\lambda}}(\ba, \bb) = \min_{\bpi \geq 0} \quad \langle \bC, \bpi\rangle
    + \lambda_1 D_\varphi(\bpi \mathds{1}_m, \ba) +  \lambda_2 D_\varphi(\bpi^\top
    \mathds{1}_n, \bb)+\lambda_{\text{reg}}D_\varphi(\bpi,\ba\bb^\top).
    \label{eq:UOTdiv_reg}
    \end{eqnarray}
As it turns out, this general problem involving different regularization weights $(\lambda_1, \lambda_2, \lambda_{\text{reg}})$ can easily be addressed in our framework as well using two simple tricks. The first one consists of absorbing the regularization weights into the divergences. Indeed, many divergences are homogeneous, i.e., satisfy a relation of the form $\lambda D_{\varphi}(\mathbf{x}|\mathbf{y})= D_{\varphi}(\lambda^{\alpha} \mathbf{x}|\lambda^{\alpha} \mathbf{y})$ where $\alpha$ is divergence-specific. This holds in particular for the KL divergence ($\alpha=1$) and the squared $\ell_{2}$ norm ($\alpha=1/2$). The second one consists of complementing $\bH$ and $\by$ with suitable terms to account for the regularization term. In the end, we may re-write Eq.~\eqref{eq:UOTdiv_reg} into Eq.~\eqref{eq:linreg} with $\lambda =1$, $\bH=[\lambda_1^{\alpha}\bH_r^\top,\lambda_2^{\alpha}\bH_c^\top,\lambda^{\alpha}_{\text{reg}}\mathbf{I}]^\top$ and $\by^\top=[\lambda^{\alpha}_1\ba^\top,\lambda^{\alpha}_2\bb^\top,\lambda^{\alpha}_{\text{reg}}\text{vec}(\ba\bb^\top)^\top]$. In particular, we obtain the following multiplicative update in the case of entropy-regularized KL-penalized UOT: 
\begin{equation}
	\bpi^{(k+1)}= \text{diag}\left(\frac{\ba}{\bpi^{(k)} \mathds{1}_m}\right)^{\frac{\lambda_1}{\lambda_{\text{all}}}}
    \left(
        \left( \bpi^{(k)}\right)^{\frac{\lambda_1+\lambda_2}{\lambda_{\text{all}}}} 
        \odot \bK
    \right) 
    \text{diag}\left(\frac{\bb}{{\bpi^{(k)\top}}\mathds{1}_n}\right)^{\frac{\lambda_2}{\lambda_{\text{all}}}}
	\label{eq:mult_update_kl_entr_reg}
\end{equation}
where $\bK=
\left(\ba\bb^\top\right)^{\frac{\lambda_{\text{reg}}}{\lambda_{\text{all}}}}
\odot \text{exp}\left(-\frac{1}{\lambda_{\text{all}}}\bC\right)$ and $\lambda_{\text{all}}=\lambda_1+\lambda_2+\lambda_{\text{reg}}$. {This multiplicative update is slightly more complex than the Sinkhorn algorithms of \citet{frogner2015learning,chizat2018scaling} and as such, it might have limited practical interest but is conceptually interesting and novel. 
Note that balanced UOT as of Eq.~\eqref{eq:UOTdiv} is simply obtained with $\lambda_{\text{reg}}=0$.} 

\paragraph{Non-linear UOT.} Finally, we discuss how our proposed reformulation of UOT can accommodate non-linear variants in which the linear term  $\langle \bC,\bpi \rangle$ is replaced by a sparsity/robustness-promoting term, leading to problems of the form
\begin{equation}
    \text{NLUOT}^{\boldsymbol{\lambda}}(\ba, \bb) = \min_{\bpi \geq 0} \quad \sum_{i,j} g(C_{i,j}T_{i,j})
+ \lambda_1 D_\varphi(\bpi \mathds{1}_m, \ba) +  \lambda_2 D_\varphi(\bpi^\top
\mathds{1}_n,\bb)
    \label{eq:nonlinearot}
\end{equation}
where $g(\cdot)$ is a usually concave function, see, e.g., \citep{candes2008enhancing,gasso2009recovering}. Our MM setting can readily accommodate such a formulation by majorizing the concave terms by their tangent. The non-linearity may improve robustness w.r.t outliers and better model realistic OT problems. For instance, in real life, the costs of transporting some goods  between two places can be nonlinear due to economies of scale.

\paragraph{Broad and potential negative societal impact.} {The contributions in this paper are methodological and focus on a  reformulation
of a fundamental OT problem and adapting existing algorithms to solve it. In
this sense, we bring more efficient solvers that run on GPU but this computational advantage can be
counterbalanced by the possibility that it brings to be applied on larger datasets. The application of OT in domain adaptation has shown that it can be used
to infer labels on samples/individuals when no labels are available, suggesting a capacity for violating user privacy. A
potential application of UOT is the case where two datasets of users acquired
by different methods contain some shared users. UOT can be used here
to find correspondences between the users in the two datasets and also identify
unique users in each dataset (those that do not receive mass). }





%
%
%

\section{Conclusion}
{In this paper, we reformulate the UOT problem as a non-negative penalized
linear regression, allowing us to propose two new classes of algorithms. We
first derive multiplicative algorithms for both KL and $\ell_2$-penalized UOT, providing numerical solutions that are fast and easy to
implement. 
For the specific case of $\ell_2$-penalized UOT, we provide the first
regularization path algorithm that computes the whole set of solutions for
\emph{all} the regularization parameter values. We finally build on the extensive
literature in inverse problem and NMF to draw some fruitful perspectives on even more efficient algorithmic solutions
or the definition of new OT problems. }
\begin{ack}
This work is partially funded by the French National Research Agency (ANR;
grants OATMIL ANR-17-CE23-0012, RAIMO ANR-20-CHIA-0021-01, MULTISCALE
ANR-18-CE23-0022-01, E4C ANR-18-EUR-0006-02, 3IA C\^ote d'Azur ANR-19-P3IA-0002,
3IA ANITI ANR-19-PI3A-0004) and the European Research Council (ERC; grant
FACTORY-CoG-6681839). 
Furthermore, this research was produced within the framework of Energy4Climate
Interdisciplinary Center (E4C) of IP Paris and Ecole des Ponts ParisTech.
This action benefited from the
support of the Chair ``Challenging Technology for Responsible Energy'' led by l'X
- Ecole Polytechnique and the Fondation de l’Ecole Polytechnique, sponsored by
TOTAL. 
\end{ack}

\bibliography{regpath}
\bibliographystyle{chicago}

\newpage
\appendix

\section{Supplementary material}
\label{sec:appendix}

\subsection{Design of $\bH$, $\bH_{r}$ and $\bH_{c}$}
\label{sec:design}


In this section, we detail how we build the design matrix $\bH$ in problem
\eqref{eq:linreg}. {By setting $\lambda = \lambda_1 = \lambda_2$},  Eq.~\eqref{eq:UOTdiv} can be reformulated
as
\begin{eqnarray}
    \text{UOT}^\lambda(\ba, \bb) = \min_{\bpi \geq 0} \quad \langle \bC, \bpi\rangle
    + \lambda D_\varphi\left(\begin{bmatrix}\bpi \mathds{1}_m \\ \bpi^\top
        \mathds{1}_n\end{bmatrix},\begin{bmatrix}\ba \\ \bb\end{bmatrix} \right) 
    \label{eq:UOTdiv2}
    \end{eqnarray}
because the divergence $D_\varphi$ is separable. Note that both $\bpi \mathds{1}_m$ and $\bpi^\top
\mathds{1}_n$ are linear operations. It means that we can vectorize the matrix $\bp=\text{vec}(\bpi)=[T_{1,1}, T_{1,2}, \dots
T_{n,m-1}, T_{n,m}]^\top$ such that:

\begin{equation}
    \begin{bmatrix}\bpi \mathds{1}_m \\ \bpi^\top
        \mathds{1}_n\end{bmatrix}=
    \bH\bp\qquad \text{where}\qquad \bH=\begin{bmatrix}
        \bH_r\\
        \bH_c
        \end{bmatrix}.
    \label{eq:H_mul_p}
    \end{equation}
The matrix $\bH_r \in \mathbb{R}_{n\times nm}$ that performs the sum over the rows of $\bpi$ is given by
\begin{equation}
    \bH_r=\begin{bmatrix}
        1 & \dots & 1 & 0 & \dots & 0& \dots & 0 & \dots & 0\\
        0 & \dots & 0 & 1 & \dots & 1& \dots & 0 & \dots & 0\\
        \hdots & \hdots & \hdots & \hdots & \hdots & \hdots & \hdots & \hdots &\hdots & \hdots\\
        0 & \dots & 0 & 0 & \dots & 0& \dots & 1 & \dots & 1\\

        \end{bmatrix}
    \label{eq:matrix_Hl}
    \end{equation}
and can be implemented in Python with $\bH_r=$ \texttt{np.repeat(np.eye(n),m)}. In a  similar fashion, the matrix
$\bH_c$ that performs the sum across columns of $\bpi$ is a $m\times nm$ array defined as
\begin{equation}
    \bH_c= \begin{bmatrix} 
       \bI_m & \bI_m & \dots & \bI_m
        \end{bmatrix} 
    \label{eq:matrix_Hc}
    \end{equation}
and can be implemented in Python using
    $\bH_c=$ \texttt{np.tile(np.eye(m),n)}.

%
%

\paragraph{Useful identities.} From the previous definitions, we have that
\begin{equation}
\bH^\top\by= 
\bH_r^\top \ba + \bH_c^\top \bb
    = \text{vec}(\ba \mathds{1}_m^\top+ \mathds{1}_n \bb^\top)= \begin{bmatrix}
        a_1+b_1\\
        a_1+b_2\\
        \hdots\\
        a_n+b_{m-1}\\
         a_n+b_{m}
        \end{bmatrix}.
    \label{eq:H_top_y}
    \end{equation}

{We have $\bH^\top\bH = \bH_r^\top\bH_r + \bH_c^\top\bH_c$, of size $nm\times nm$. $\bH_r^\top\bH_r$ is a block-diagonal matrix with $n$ blocks of size $m \times m$ filled with ones. $\bH_r^\top\bH_r$ can be implemented in Python with \texttt{np.tile(np.eye(m),(n,m))}. $\bH_c^\top\bH_c$ is a block matrix with blocks of $ \bI_m$, and can be implemented in Python with \texttt{np.tile(np.eye(m), (n,n))}. Multiplying $\bH^\top\bH$ by a vector, e.g. $\bp$, results in $\bH^\top\bH\bp = \text{vec}(\bpi \mathds{1}_m \mathds{1}_m^\top+  \mathds{1}_n \mathds{1}_n^\top\bpi )$.}

\subsection{Details of MM algorithms}
\label{sec:app_mult}

The objective function $F_{\lambda}(\bp)$ defined by Eq.~\eqref{eq:linreg} can be re-written as:
\begin{equation}
F_{\lambda}(\bp) = \sum_{i} \varphi(\sum_{j} H_{i,j} t_{j}) + \sum_{j} \left[ \frac{c_{j}}{\lambda} - \sum_{i} H_{i,j} \varphi'(y_{i})\right] t_{j}.
\end{equation}
Applying Jensen inequality to the first term like explained in Section~\ref{sec:mmforuot} directly leads to the expression of $G_{\lambda}(\bp,\tilde{\bp})$ given by Eq.~\eqref{eq:auxiliary_func}. The auxiliary function is separable and convex. Given $\tilde{\bp} = {\bp}^{(k)}$, the next iterate $\bp^{(k+1)}$ can be computed by cancelling the partial derivative $\nabla_{t_q} G_\lambda(\bp,\bp^{(k)})$, $q=1,\ldots,nm$, or setting $t_q$ to zero if the solution is negative in order to satisfy the non-negative constraint (note that this is not a heuristic but what the KKT conditions dictate). Cancelling the partial derivative w.r.t. $t_{q}$ is equivalent to solving
\begin{eqnarray}
\sum_{i} H_{i,q} \varphi'\left( \frac{t_q}{t^{(k)}_q} [\bH \bp^{(k)}]_i \right) = \sum_{i} H_{i,q} \varphi'(y_{i}) -  \frac{c_{q}}{\lambda}
\label{eq:opt_cdt_aux_fct}
\end{eqnarray}
w.r.t. $t_{q}$. We address this univariate problem for the $\ell_2$ and KL-penalties next. 

\paragraph{Squared $\ell_2$ penalty.} In that case we have $\varphi (x) = \frac{x^2}{2} $,  $ \varphi'(x) = x$ and we obtain
\begin{equation}
	t_q^{(k+1)} =  t^{(k)}_q \, \frac{\max \left(0, [\bH^\top \by]_q - \frac{1}{\lambda} c_{q}\right)}{[\bH^\top  \bH \bp^{(k)} ]_q}.
	\label{eq:general_Quad_multiplicative_update}
\end{equation}
Recall that $\bp$ is a vector form of the OT plan $\bpi$, and assume that $t_q$ corresponds to the entry $T_{i,j}$.  $\bH^\top \by$ is a $nm$-dimensional vector with elements $a_i + b_j$, see Eq. \eqref{eq:H_top_y}. Furthermore, we have $\bH \bp =  \begin{bmatrix}\bpi \mathds{1}_m \\ \bpi^\top
	\mathds{1}_n\end{bmatrix} $ thanks to Eq.  \eqref{eq:H_mul_p}. Therefore, we can establish the following update in $T_{i,j}$
\begin{equation}
	T_{i,j}^{(k+1)} =  T_{i,j}^{(k)} \, \frac{\max \left(0, a_i + b_j - \frac{1}{\lambda} c_{i,j}\right)}{[\bpi^{(k)} \mathds{1}_m ]_i + [{\bpi^{(k)}}^\top
		\mathds{1}_n]_j}
\label{eq:MML2}
\end{equation}
with matrix form given by Eq.~\eqref{eq:mult_update_quadratic_divergence_l2}. 
 
 \paragraph{KL penalty.} In this case we have $\varphi (x) = x \log x - x $, $\varphi' (x) = \log x$ and we obtain
 \begin{align}
 t^{(k+1)}_{q} &= t_{q}^{(k)} \exp \left( \frac{1}{\sum_{q} H_{i,q}} \left (\sum_{i} H_{i,q} \log \frac{y_{i}}{[\bH \bp^{(k)}]_{i}}  -\frac{c_{q}}{\lambda} \right) \right) \\
  	 &= t^{(k)}_q \, \exp\left(\frac{\left[\bH^\top \log(\by) - \bH^\top \log\left(  \bH \bp^{(k)} \right)\right]_q - \frac{1}{\lambda} c_{q}}{\left[\bH^\top \mathds{1} \right]_q}\right). \label{eq:mmupdate_kl_H}
  \end{align}

Using the results of Section~\ref{sec:design} like in the $\ell_{2}$ case, we obtain the following update
 \begin{eqnarray*}
T_{i,j}^{(k+1)} & =  & \left(\frac{a_i}{[\bpi^{(k)} \mathds{1}_m]_i}\right)^{1/2} \, T_{i,j}^{(k)} \exp\left(- \frac{c_{i,j}}{2 \lambda}\right)  \, \left(\frac{b_j}{[{\bpi^{(k)}}^\top \mathds{1}_n]_j}\right)^{1/2}
\label{eq:MMKL}
\end{eqnarray*}
with matrix form given by Eq.~\eqref{eq:mult_update_kl}.

\paragraph{Alternative multiplicative update for the $\ell_2$-penalty.}
{Another possible approach is to use a quadratic majorization of the linear term $\bc^{\top} \bp$ to bypass the thresholding operation like in \citep{hoyer2002non,yang11}, leading to: 
\begin{equation}
	\bpi^{(k+1)}= \bpi^{(k)} \odot \frac{ \ba \mathds{1}_m^\top + \mathds{1}_n \bb^\top   }{\bpi^{(k)} \mathbf{O}_m + \mathbf{O}_n \bpi^{(k)} +\frac{1}{2\lambda} \bC} \quad \text{with} \quad \mathbf{O}_\ell = \mathds{1}_\ell \mathds{1}_\ell^\top.
	\label{eq:mult_update_quadratic_divergence_l22}
\end{equation}
However we found update~\eqref{eq:mult_update_quadratic_divergence_l2} more useful in our case, thanks to the thresholding operation that locates true zeros from start.}

\paragraph{{Alternative derivation of MM algorithms.}} The reformulation of UOT as a non-negative penalized linear regression problem comes very handy because it offers a novel interpretation of UOT and the possibility of using some of the many existing algorithms for the latter problem, such as LARS-based algorithm for path computation. However, we want to point out that we may also derive MM algorithms directly from~Eq.~\eqref{eq:UOTdiv}. Let us write
\bal{
F_{\lambda}(\bpi) &= \langle \bC, \bpi\rangle
+ \lambda_1 D_\varphi(\bpi \mathds{1}_m, \ba) +  \lambda_2 D_\varphi(\bpi^\top
\mathds{1}_n, \bb) \\
&= \sum_{ij} C_{i,j} T_{i,j} + \lambda_{1} \sum_{i} d_{\varphi}(\sum\nolimits_{j} T_{i,j} , a_{i}) + \lambda_{2} \sum_{j} d_{\varphi}(\sum\nolimits_{i} T_{i,j} , b_{j}) \label{eq:altF}
}
(Note that we have $F_{\lambda}(\bpi) = F_{\lambda}(\bp)$, slightly abusing notations). Let $\tilde{\bpi}$ be a current estimate of $\bpi$. We wish to compute an auxiliary function $G_{\lambda}(\bpi,\tilde{\bpi})$ for $F_{\lambda}(\bpi)$.
Let us denote
\bal{
\tilde{a}_{i} &= \sum_{j} \tilde{T}_{i,j} \quad \text{(the $i^{th}$ approximate row marginal)} \\
\tilde{b}_{j} &= \sum_{i} \tilde{T}_{i,j} \quad \text{(the $j^{th}$ approximate column marginal)} \\
\tilde{\alpha}_{i,j} & = \frac{\tilde{T}_{i,j}}{\tilde{a}_{i}} \quad \ \ \quad \text{such that } \sum_{j} \tilde{\alpha}_{i,j} =1 \\
\tilde{\beta}_{i,j} & = \frac{\tilde{T}_{i,j}}{\tilde{b}_{j}} \quad \ \ \quad \text{such that } \sum_{i} \tilde{\beta}_{i,j} =1
}
By convexity of $d_{\varphi}(x,y)$ w.r.t $x$, we have
\bal{
d_{\varphi}\left(\sum_{j} T_{i,j} , a_{i}\right) & \le \sum_{j} \tilde{\alpha}_{i,j} \, d_{\varphi}\left( \frac{T_{i,j}}{\tilde{\alpha}_{i,j}} ,a_{i} \right), \\
d_{\varphi}\left(\sum_{i} T_{i,j} , b_{j}\right) & \le \sum_{i} \tilde{\beta}_{i,j} \, d_{\varphi}\left( \frac{T_{i,j}}{\tilde{\beta}_{i,j}} , b_{j}\right). 
}
The inequalities are tight when $\tilde{\bpi} = \bpi$. Plugging the latter inequalities into Eq.~\eqref{eq:altF}, we obtain the following auxiliary function:
\bal{
G_{\lambda}(\bpi|\tilde{\bpi}) = \sum_{ij} \left[ C_{i,j} T_{i,j} + \lambda_{1} \tilde{\alpha}_{i,j} \, d_{\varphi}\left( \frac{T_{i,j}}{\tilde{\alpha}_{i,j}} ,a_{i} \right) + \lambda_{2} \tilde{\beta}_{i,j} \, d_{\varphi}\left( \frac{T_{i,j}}{\tilde{\beta}_{i,j}} , b_{j}\right) \right].
}
$G_{\lambda}(\bpi|\tilde{\bpi})$ is essentially the matrix form of $G_{\lambda}(\bp|\tilde{\bp})$, with partial derivative given by: 
\bal{
\nabla_{T_{i,j}} G_{\lambda}(\bpi|\tilde{\bpi}) = C_{i,j} + \lambda_{1} d_{\varphi}'\left(\tilde{a}_{i} \frac{T_{i,j}}{\tilde{T}_{i,j}}  , a_{i}\right) + \lambda_{2} d_{\varphi}'\left(\tilde{b}_{j} \frac{T_{i,j}}{\tilde{T}_{i,j}}  , b_{j}\right).
}
Using $d'_{\varphi}(x,y) = \varphi'(x) - \varphi'(y)$ and either $\varphi'(x)=x$ ($\ell_{2}$-penalized UOT) or $\varphi'(x)=\log x$ (KL-penalized UOT), we easily retrieve Eq.~\eqref{eq:MML2} and Eq.~\eqref{eq:MMKL} when $\lambda_{1}=\lambda_{2}$, or Eq.~\eqref{eq:mult_update_kl_entr_reg} in the general case (with here $\lambda_{\text{reg}} = 0$).

\subsection{Details of the UOT path computation}
\label{sec:matrices_reg_path}

\paragraph{Matrices and vectors on the active set $\setA$.} Recall that $\bm_{\setA}$,  $\bc_{\setA}$ and  $\bp_{\setA}$ are sub-vectors of $\bm$,  $\bc$ and  $\bpi$ corresponding to indices in $\setA$. $\bH_{\setA}$ is a matrix of dimension $(|i| + |j|)\times |\setA|$, where $|i|$ and $|j|$ are respectively the number of distinct rows $i$ and columns $j$ that belong to the transport plan for a given active set $\setA$. $\bH_{\setA}$ is built by keeping only the rows of $\bH_r$ such that the element $i$ is present in the active set (the latter matrix being denoted $[\bH_{r}]_\setA$), the rows of $\bH_c$ such that the element $j$ is present in the active set (denoted $[\bH_{c}]_\setA$), and keeping the columns such that element 
$(i,j)\in \setA$ ({up to vectorization}).

\paragraph{Update $(\bH_{\setA}^\top\bH_{\setA})^{-1}$ from $(\bH_{\setA_k}^\top\bH_{\setA_k})^{-1}$ using the Schur complement.}

Algorithm \ref{alg:regpath} involves the computation, 
{at each iteration, of the inverse  matrix} $(\bH_{\setA}^\top\bH_{\setA})^{-1}$. The computational burden can be alleviated by using the Schur complement of the matrix in order to compute $(\bH_{\setA}^\top\bH_{\setA})^{-1}$ from its value at the previous iteration $(\bH_{\setA_k}^\top\bH_{\setA_k})^{-1}$. Let us denote $\bB_{\setA} = (\bH_{\setA}^\top\bH_{\setA})$ and $\bB_{\setA_k} = (\bH_{\setA_k}^\top\bH_{\setA_k})$. Two cases may arise:
\begin{itemize}
\item One component $q$ is added to the active set $\setA = \setA_{k+1}= \setA_{k} \cup q$. In that case, we have:
\begin{equation}
\bB^{-1}_{\setA} = \begin{bmatrix}
\bB^{-1}_{\setA_{k}} + \bB^{-1}_{\setA_k}b_{\setA, q}S^{-1}b_{q,\setA} \bB^{-1}_{\setA_{k}} &  -\bB^{-1}_{\setA_{k}} b_{\setA, q}S^{-1}\\
 - S^{-1}b_{{q},\setA} \bB^{-1}_{\setA}  & S^{-1}
\end{bmatrix}
\end{equation}
where $b_{q,\setA}$ is the last column of matrix $\bB_{\setA}$, $b_{\setA,q}$ its  last row and $S = 2 - b_{q,\setA}{^\top}\bB^{-1}_{\setA_{k}}b_{\setA,q}$ is a scalar. 
\item One component $q$ is removed from the active set $\setA = \setA_{k} \backslash  q$. In that case, we get:
\begin{equation}
\bB^{-1}_{\setA} = \bB^{-1}_{\setA_{k}\backslash  q} - \dfrac{b^{-1}_{\setA\backslash q, q}b^{-1}_{q, \setA\backslash q}}{b^{-1}_{q,q}} \end{equation}
with $\bB^{-1}_{\setA\backslash q}$ being the matrix $\bB^{-1}_{\setA}$ deprived from its row and column corresponding to the component $q$. The vector  $b^{-1}_{\setA_{k}\backslash  q,q}$ represents the column of the $\bB^{-1}_{\setA}$ matrix corresponding to element $i$ while $b^{-1}_{q, \setA\backslash  q}$ stands for the  corresponding row. Finally $b^{-1}_{q,q}$ is the component of $\bB^{-1}_{\setA}$ corresponding to the component  $q$. 
\end{itemize}

\subsection{Details of the regularization path formulation for semi-relaxed UOT}
\label{sec:reg_path_semi}

\paragraph{Semi-relaxed $\ell_2$-penalized UOT.} 
We start by recalling the formulation of the semi-relaxed $\ell_2$-penalized UOT problem:
$$
	\text{SROT}^{\boldsymbol{\lambda}}(\ba, \bb) = \min_{\bpi \geq 0, \bH_c \bp = \bb} \,\, \langle \bC, \bpi\rangle
	+ \lambda \|\bpi \mathds{1}_m -  \ba\|^2.
	$$ 
	
From Eq.~\eqref{eq:lagrang_SROT}, the corresponding Lagrangian writes:
\begin{equation}
	{L}_\lambda(\bp, \bgamma) = \frac{1}{\lambda}\bc^\top\bp + {\frac{1}{2}}(\bH_r\bp-\ba)^\top(\bH_r\bp-\ba)+ (\bH_c\bp-\bb)^\top\bu-\bgamma^\top\bp
\end{equation}
with $\bu \in \mathbb{R}^m$ the Lagrange parameters associated to the $m$ equality constraints and $\bgamma \geq 0$ the Lagrange parameters related to the non-negativity constraints. We recall the KKT optimality conditions, which state that i) $\nabla_{\bp} {L}_{\lambda}=\frac{1}{\lambda}\bc + \bH_r^\top\bH_r\bp-\bH_r^\top\ba +\bH_c^\top\bu- \bgamma = 0$ (stationary condition), ii) $\bgamma  \odot \bp = 0$ (complementary condition), and iii) $\bgamma \geq 0$  and $\bH_c\bp-\bb=\boldsymbol{0}$ (feasibility) from which we may derive  the path computation. We recall that $\odot$ stands for point-wise multiplication.

\paragraph{Piecewise linearity of the path.} Let us suppose that, at step $k$, we know the current active set $\setA = \mathcal{A}_k$ and we look for $\bp_{\setA}^{\lambda}$ and $\bu^{\lambda}$. 
Because of the complementary condition, we have $\bgamma_{\setA} =
\boldsymbol{0}$. Hence the stationnarity condition on the active set can be
rewritten as, with $\lambda = \lambda_k+\epsilon$  and $\epsilon$ small enough 
\begin{eqnarray}
\begin{cases}
[\bH_{r}^\top]_{\setA}[\bH_{r}]_{\setA}\bp^{\lambda}_{\setA} +  [\bH_{c}^\top]_{\setA}\bu^{\lambda} &= [\bH_{r}^\top]_{\setA}\ba_{\setA} - \frac{1}{\lambda}\bc_{\setA} \\ 
[\bH_{c}]_{\setA}\bp^{\lambda}_{\setA}&=\bb_{\setA}
\end{cases}
\end{eqnarray}
or equivalently, at each iteration, the following linear system  should be solved:
\begin{equation}
\label{ls_matrix}
	\underbrace{\begin{pmatrix}
	 [\bH_{r}^\top]_{\setA}[\bH_{r}]_{\setA} &[\bH_{c}^\top]_{\setA}  \\
		[\bH_{c}]_{\setA}& \boldsymbol{0}
	\end{pmatrix}}_{{\bK}_{\setA}}
\begin{pmatrix}
	\bp^{\lambda}_{\setA} \\
	\bu^{\lambda}
\end{pmatrix}
= 
-\frac{1}{\lambda}\underbrace{\begin{pmatrix}
	\bc_{\setA} \\
	\boldsymbol{0}
\end{pmatrix}}_{\boldsymbol{\gamma}_{\setA}}  +
\underbrace{
\begin{pmatrix}
	 [\bH_{r}^\top]_{\setA}\ba_{\setA} \\
	\bb
\end{pmatrix}}_{\boldsymbol{\beta}_{\setA}}.
\end{equation}
We then have
\begin{equation}
\begin{pmatrix}
	\bp^{\lambda}_{\setA} \\
	\bu^{\lambda}
\end{pmatrix}
= 
-\frac{1}{\lambda}\bK_{\setA}^{-1}\boldsymbol{\gamma}_{\setA}  + \bK_{\setA}^{-1}
{\boldsymbol{\beta}_{\setA}}.
\end{equation}
We now denote $\tilde{\bc}_{\setA} = \bK_{\setA}^{-1}\boldsymbol{\gamma}_{\setA}$ and its sub-vectors $\tilde{\bc}_{\setA}^{\ba}$ and  $\tilde{\bc}_{\setA}^{\bb}$ that respectively contains the $|\setA|$ first rows and $m$ last rows of $\tilde{\bc}_{\setA}$. We also denote $\tilde{\bm}_{\setA} = \bK_{\setA}^{-1}{\boldsymbol{\beta}_{\setA}}$ and its sub-vectors $\tilde{\bm}_{\setA}^{\ba}$ and $\tilde{\bm}_{\setA}^{\bb}$ in the same fashion. We then have 
\begin{eqnarray}
\begin{cases}
\bp^{\lambda}_{\setA} = -\frac{1}{\lambda}\tilde{\bc}_{\setA}^{\ba}  + \tilde{\bm}_{\setA}^{\ba} \\
\bu^{\lambda} = -\frac{1}{\lambda}\tilde{\bc}_{\setA}^{\bb}  + \tilde{\bm}_{\setA}^{\bb} \\ 
\end{cases}
\label{eq:solsemi}
\end{eqnarray}
We again notice the piecewise linearity (as a function of ${1}/{\lambda}$) of the path when the active set $\mathcal{A}$ is fixed.

\paragraph{Computation of $\lambda^{k+1}$ given $\lambda^k$.}
Given a current solution at iteration $k$ $(\lambda_k, \bp^{\lambda_k})$, we increase the $\epsilon$ value in $\lambda = \lambda_k+\epsilon$ until one of the following case arises. 

$\bullet$ \textbf{Inside the active set,} the positivity constraint on $\bp_{\setA}^{\lambda}$ may be violated, corresponding to the case 
\begin{eqnarray}
\tilde{\bm}_{\setA}^{\ba} = \frac{1}{\lambda} \tilde{\bc}_{\setA} ^{\ba}
 & \Rightarrow &
  \lambda_r = \min_{>\lambda_k} \left( \frac{\tilde{\bc}_{\setA}^{\ba}}{ \tilde{\bm}_{\setA}^{\ba}  } \right)
\end{eqnarray}
where $\min_{>\lambda_k}$ denotes the smallest value in $\frac{\tilde{\bc}_{\setA}^{\ba}}{ \tilde{\bm}_{\setA}^{\ba}  }$ greater that $\lambda_k$.

$\bullet$ \textbf{Outside the active set,} the positivity constraint of the KKT may be violated. The stationarity condition outside the active set $\bar{\setA}$ can be rewritten, by injecting the solution of Eq. \eqref{eq:solsemi}:

\begin{equation}
	 \frac{1}{\lambda}\bc_{\bar{\setA}}+  \big[\bH_r^\top(\bH_r(-\frac{1}{\lambda}\tilde{\bc}^{\ba}  + \tilde{\bm}^{\ba})-\ba)\big]_{\bar{\setA}}+ \big[\bH_c^\top( -\frac{1}{\lambda}\tilde{\bc}^{\bb}  + \tilde{\bm}^{\bb})\big]_{\bar{\setA}}-\bgamma_{\bar{\setA}} = 0
\end{equation}

\begin{eqnarray}
\frac{1}{\lambda}\bc_{\bar{\setA}}+  \Big[\bH^\top \bH\big(\tilde{\bm} + \frac{1}{\lambda}\tilde{\bc}\big)\Big]_{\bar{\setA}}-    \bm_{\bar{\setA}} = \bgamma_{\bar{\setA}}
& \Rightarrow
&
\lambda_a =  \min_{>\lambda_k}\left(\dfrac{\bc_{\bar{\setA}} - \big[\bH^\top\bH\tilde{\bc}\big]_{\bar{\setA}} }{\bm_{\bar{\setA}}- \big[\bH^\top\bH\tilde{\bm}\big]_{\bar{\setA}} }\right) 
\end{eqnarray}
The active set changes only if there exists a component $i$ outside the current active set such that $\gamma_{i}\geq 0$.  Hence we write:
\begin{eqnarray}
	 \frac{1}{\lambda}\bc_{\bar{\setA}} - \frac{1}{\lambda}\big[\bH_r^\top\bH_r\tilde{\bc}^{\ba} + \bH_c^\top\tilde{\bc}^{\bb}\big]_{\bar{\setA}} +  \big[\bH_r^\top\bH_r \tilde{\bm}^{\ba}-\bH_r^\top\ba + \bH_c^\top\tilde{\bm}^{\bb}\big]_{\bar{\setA}}\geq 0 \\
{\lambda_a} = \min_{>\lambda_k}\dfrac{\bc_{\bar{\setA}} - \big[\bH_r^\top\bH_r\tilde{\bc}^{\ba} + \bH_c^\top\tilde{\bc}^{\bb}\big]_{\bar{\setA}}}{\big[2\bH_r^\top\ba -\bH_r^\top\bH_r \tilde{\bm}^{\ba} - \bH_c^\top\tilde{\bm}^{\bb}\big]_{\bar{\setA}}}
\end{eqnarray}
Note that this last equation is very similar to the one we obtain for $\ell_2$-penalized UOT, except that vectors $\tilde{\bm}$ and $\tilde{\bc}$ are split in 2 parts, depending on if we consider the rows (that can be unbalanced) or the columns (that should strictly respect the marginal constraint). Also note that the Schur complement applies to the update of $\bK_{\setA}^{-1}$ in order to decrease the computational burden.



\subsection{IPOT is a MM algorithm}
\label{sec:prox_mm}

{Herein we discuss the relation between the Inexact
Proximal Point OT (IPOT) algorithm of \citet{xie2020fast} and MM. First note that IPOT aims at 
the balanced OT problem \eqref{eq:exactot}. This is equivalent to solving
\begin{equation}
    \label{eq:primalOTDC}
 \min_{\bpi \geq 0, \bpi \mathds{1}_m = \ba, \bpi^\top \mathds{1}_n = \bb} \langle \bC, \bpi\rangle\quad + \lambda \sum_{i,j} T_{i,j}\log(T_{i,j})-\lambda \sum_{i,j} T_{i,j}\log(T_{i,j})
  \end{equation}
where one adds and removes the entropy regularization of $\bpi$. A simple algorithm can be devised by upper bounding the concave term by its tangent at $\bpi^{(k)}$ leading to the new problem
%
%
%
  \begin{equation}
    \label{eq:dca}
    \min_{\bpi \geq 0, \bpi \mathds{1}_m = \ba, \bpi^\top \mathds{1}_n = \bb}  \langle \bpi,\bC \rangle +\lambda \sum_{i,j} T_{i,j}\log(T_{i,j})-\lambda \langle \bpi,\log\left(\bpi^{(k)}\right)+1 \rangle 
  \end{equation}
  where the log is taken component-wise. Note that the constant $1$ in the
  scalar product can be removed since $\sum _{i,j}T_{i,j}$ is constant and does not
  influence the solution. Problem~\eqref{eq:dca} can be solved using classical Sinkhorn iterations with 
   a cost matrix $\tilde\bC=\bC-\lambda\log(\bpi^{(k)})$. This corresponds to using the kernel matrix
  \begin{equation}
    \tilde\bK=\exp\left(-\frac{1}{\lambda}(\bC-\lambda\log\left(\bpi^{(k)}\right)\right)=\exp\left(-\frac{1}{\lambda}\bC\right)\odot\bpi^{(k)},
      \label{eq:kernelipot}
  \end{equation}
as presented in \cite[Algorithm 1]{xie2020fast}. Hence IPOT can be interpreted as MM. Note that the point-wise product between a kernel matrix and the estimate $\bpi^{(k)}$ appears also in our multiplicative updates~\eqref{eq:mult_update_kl} and \eqref{eq:mult_update_kl_entr_reg}
with however a different scaling parameter. 

}

\subsection{Details about the experiments}
\label{sec:details_expes}
We run the experiments on a Mac mini 2020 personal computer, with M1 chip and 16GB of RAM. All the experiments can be re-run thanks to the paper companion code. We compare the following algorithms provided by the following solvers:
\begin{itemize}
\item the ``L-BFGS-B'' method of SciPy, in which we provide the function to minimize and its associated Jacobian (either for KL or  $\ell_2$-penalized UOT),
\item the Lasso algorithms Celer and of Scikit-learn,
\item the regularization path algorithm introduced in the paper,
\item the multiplicative updates introduced in the paper.
\end{itemize}
We use the same stopping criteria for all the algorithms (not to mention the regularization path algorithm that provides an exact solution), except for $\ell_2$-penalized UOT that necessitates a smaller tolerance to converge to the correct values, especially for large values of $\lambda$. 

Regarding Figure \ref{fig:timings}, we draw  5 realizations of two random 2 Gaussian samples of 10-dimensional $n=m$ points with different means and variances. 

\end{document}